\documentclass{article}

\usepackage{amsmath, amsfonts, amsthm, amssymb}

\usepackage{enumitem}

\usepackage{times}
\usepackage{graphicx} 
\usepackage{subfigure} 

\usepackage{natbib}

\usepackage{algorithm}
\usepackage{algorithmic}

\PassOptionsToPackage{hyphens}{url}
\usepackage{hyperref}

\allowdisplaybreaks

\newtheorem{theorem}{Theorem}
\newtheorem{lemma}[theorem]{Lemma}

\newtheorem{corollary}[theorem]{Corollary}
\newtheorem{definition}[theorem]{Definition}
\newtheorem{remark}[theorem]{Remark}

\newcommand{\acro}[1]{\textsc{\MakeLowercase{#1}}}

\newcommand{\inv}{^{-1}}                 
\newcommand{\sminus}{\backslash}         
\newcommand{\N}{\mathbb{N}}              
\newcommand{\R}{\mathbb{R}}              
\newcommand{\e}{\varepsilon}             
\renewcommand{\d}{\delta}                
\newcommand{\X}{\mathcal{X}}             
\newcommand{\XX}{\mathbb{X}}             
\newcommand{\E}{\mathop{\mathbb{E}}}     
\newcommand{\Var}{\mathbb{V}}            
\newcommand{\pr}{\mathbb{P}}             

\usepackage[accepted]{icml2016}

\icmltitlerunning{Analysis of $k$-Nearest Neighbor Statistics with Application
to Entropy Estimation}

\begin{document} 

\twocolumn[
\icmltitle{Analysis of $k$-Nearest Neighbor Distances \\ with Application to
Entropy Estimation}

\icmlauthor{Shashank Singh}{sss1@andrew.cmu.edu}
\icmlauthor{Barnab\'as P\'oczos}{bapoczos@cs.cmu.edu}
\icmladdress{Carnegie Mellon University,
             5000 Forbes Ave., Pittsburgh, PA 15213 USA}

\icmlkeywords{entropy, nonparametric statistics, nearest neighbor}

\vskip 0.3in
]

\begin{abstract}
Estimating entropy and mutual information consistently is important for many
machine learning applications. The Kozachenko-Leonenko (KL) estimator
\citep{kozachenko87statistical} is a widely used nonparametric estimator for
the entropy of multivariate continuous random variables, as well as the basis
of the mutual information estimator of \citet{Kraskov04estimating}, perhaps the
most widely used estimator of mutual information in this setting. Despite the
practical importance of these estimators, major theoretical questions regarding
their finite-sample behavior remain open. This paper proves finite-sample
bounds on the bias and variance of the KL estimator, showing that it achieves
the minimax convergence rate for certain classes of smooth functions. In
proving these bounds, we analyze finite-sample behavior of $k$-nearest
neighbors ($k$-NN) distance statistics (on which the KL estimator is based). We
derive concentration inequalities for $k$-NN distances and a general
expectation bound for statistics of $k$-NN distances, which may be useful for
other analyses of $k$-NN methods.
\end{abstract} 

\section{Introduction}
\label{sec:intro}

Estimating entropy and mutual information in a consistent manner is of
importance in a number problems in machine learning. For example, entropy
estimators have applications in
goodness-of-fit testing \cite{goria05new},
parameter estimation in semi-parametric models \cite{Wolsztynski85minimum},
studying fractal random walks \cite{Alemany94fractal},
and texture classification \cite{hero02alpha,hero2002aes}.
Mutual information estimators have applications in
feature selection \cite{peng05feature},
clustering \cite{aghagolzadeh07hierarchical},
causality detection \cite{Hlavackova07causality},
optimal experimental design \cite{lewi07realtime, poczos09identification},
f\acro{MRI} data processing \cite{chai09exploring},
prediction of protein structures \cite{adami04information},
and boosting and facial expression recognition \cite{Shan05conditionalmutual}.
Both entropy estimators and mutual information estimators have been used for
independent component and subspace analysis
\cite{radical03,szabo07undercomplete_TCC, poczos05geodesic,Hulle08constrained},
as well as for image registration
\cite{kybic06incremental,hero02alpha,hero2002aes}.
For further applications, see
\cite{Leonenko-Pronzato-Savani2008}.

In this paper, we focus on the problem of estimating the Shannon entropy of
a continuous random variable given samples from its distribution. All of our
results extend to the estimation of mutual information, since the latter can be
written as a sum of entropies.
\footnote{Specifically, for random variables $X$ and $Y$,
$I(X; Y) = H(X) + H(Y) - H(X, Y)$.}
In our setting, we assume we are given $n$ IID
samples from an unknown probability measure $P$. Under nonparametric
assumptions (on the smoothness and tail behavior of $P$), our task is then to
estimate the differential Shannon entropy of $P$.

Estimators of entropy and mutual information come in many forms (as reviewed
in Section \ref{sec:related}), but one common approach is based on statistics
of $k$-nearest neighbor ($k$-NN) distances (i.e., the distance from a sample to
its $k^{th}$ nearest neighbor amongst the samples, in some metric on the
space). These nearest-neighbor estimates are largely based on initial work by
\citet{kozachenko87statistical}, who proposed an estimate for differential
Shannon entropy and showed its weak consistency. Henceforth, we refer to this
historic estimator as the `KL estimator', after its discoverers. Although there
has been much work on the problem of entropy estimation in the nearly three
decades since the KL estimator was proposed, there are still major open
questions about the finite-sample behavior of the KL estimator. The goal of
this paper is to address some of these questions in the form of finite-sample
bounds on the bias and variance of the estimator.

Specifically, our {\bf main contributions} are the following:
\begin{enumerate}
\item
We derive 
$O \left( \left( k/n \right)^{\beta/D} \right)$ bounds on the bias of the KL
estimate, where $\beta$ is a measure of the smoothness (i.e., H\"older
continuity) of the sampling density, $D$ is the intrinsic dimension of the
support of the distribution, and $n$ is the sample size.
\item
We derive $O \left( n\inv \right)$ bounds on the variance of the KL estimator.
\item
We derive concentration inequalities for $k$-NN distances, as well as general
bounds on expectations of $k$-NN distance statistics, with important special
cases:
\begin{enumerate}
\item
We bound the moments of $k$-NN distances, which play a role in analysis of many
applications of $k$-NN methods, including both the bias and variance of the KL
estimator. In particular, we significantly relax strong assumptions underlying
previous results by \citet{evans02KNNmoments}, such as compact support and
smoothness of the sampling density. Our results are also the first which apply
to negative moments (i.e., $\E \left[ X^\alpha \right]$ with $\alpha < 0$);
these are important for bounding the variance of the KL estimator.
\item
We give upper and lower bounds on the logarithms of $k$-NN distances. These are
important for bounding the variance of the KL estimator, as well as $k$-NN
estimators for divergences and mutual informations.
\end{enumerate}
\end{enumerate}

We present our results in the general setting of a set equipped with a metric,
a base measure, a probability density, and an appropriate definition of
dimension. This setting subsumes Euclidean spaces, in which $k$-NN methods have
traditionally been analyzed,
\footnote{A recent exception in the context of classification, is
\citet{chaudhuri14KNNrates} which considers general metric spaces.}
but also includes, for instance, Riemannian manifolds, and perhaps other
spaces of interest. We also strive to weaken some of the restrictive
assumptions, such as compact support and boundedness of the density, on which
most related work depends.

We anticipate that the some of the tools developed here may be used to derive
error bounds for $k$-NN estimators of mutual information, divergences
\citep{Wang-Kulkarni-Verdu2009}, their generalizations (e.g., R\'enyi and
Tsallis quantities \citep{Leonenko-Pronzato-Savani2008}), norms, and other
functionals of probability densities. We leave such bounds to future work.

\subsection*{Organization}
Section \ref{sec:related} discusses related work. Section \ref{sec:setting}
gives theoretical context and assumptions underlying our work. In Section
\ref{sec:concentration}, we prove concentration boundss for $k$-NN distances,
and we use these in Section \ref{sec:KNN_stats} to derive bounds on the
expectations of $k$-NN distance statistics. Section \ref{sec:KL_est} describes
the KL estimator, for which we prove bounds on the bias and variance in
Sections \ref{sec:bias_bound} and \ref{sec:variance_bound}, respectively.

\section{Related Work}
\label{sec:related}
Here, we review previous work on the analysis of $k$-nearest neighbor
statistics and their role in estimating information theoretic functionals, as
well as other approaches to estimating information theoretic functionals.

\subsection{The Kozachenko-Leonenko Estimator of Entropy}

In general contexts, only weak consistency of the KL estimator is known
\citep{kozachenko87statistical}.
\citet{biau15EntropyKNN} recently reviewed finite-sample results known for
the KL estimator. They show (Theorem 7.1) that, if the density $p$ has compact
support, then the variance of the KL estimator decays as $O(n\inv)$. They also
claim (Theorem 7.2) to bound the bias of the KL estimator by $O(n^{-\beta})$,
under the assumptions that $p$ is $\beta$-H\"older continuous
($\beta \in (0, 1]$), bounded away from $0$, and supported on the interval
$[0, 1]$. However, in their proof \citet{biau15EntropyKNN} neglect the
additional bias incurred at the boundaries of $[0, 1]$, where the density cannot
simultaneously be bounded away from $0$ and continuous. In fact, because the KL
estimator does not attempt to correct for boundary bias, for densities bounded
away from $0$, the estimator may suffer bias worse than $O(n^{-\beta})$.

The KL estimator is also important for its role in the mutual information
estimator proposed by \citet{Kraskov04estimating}, which we refer to as the KSG
estimator. The KSG estimator expands the mutual information as a sum of
entropies, which it estimates via the KL estimator with a particular
\emph{random} (i.e., data-dependent) choice of the nearest-neighbor parameter
$k$. The KSG estimator is perhaps the most widely used estimator for the mutual
information between continuous random variables, despite the fact that it
currently appears to have no theoretical guarantees, even asymptotically. In
fact, one of the few theoretical results, due to
\citet{gao15stronglyDependent}, concerning the KSG estimator is a negative
result: when estimating the mutual information between strongly dependent
variables, the KSG estimator tends to systematically underestimate mutual
information, due to increased boundary bias.
\footnote{To alleviate this, \citet{gao15stronglyDependent} provide a heuristic
correction based on using local PCA to estimate the support of the
distribution. \citet{gao15localGaussian} provide and prove asymptotic
unbiasedness of another estimator, based on local Gaussian density estimation,
that directly adapts to the boundary.}
Nevertheless, the widespread use of the KSG estimator motivates study of its
behavior. We hope that our analysis of the KL estimator, in terms of which the
KSG estimator can be written, will lead to a better understanding of the
latter.


\subsection{Analysis of nearest-neighbor distance statistics}

\citet{evans08SLLNforKNN} derives a law of large numbers for $k$-NN statistics
with uniformly bounded (central) kurtosis as the sample size  $n \to \infty$.
Although it is not obvious that the kurtosis of $\log$-$k$-NN distances is
uniformly bounded (indeed, each $\log$-$k$-NN distance approaches $-\infty$
almost surely), we show in Section \ref{sec:variance_bound} that this is indeed
the case, and we apply the results of \citet{evans08SLLNforKNN} to bound the
variance of the KL estimator.

\citet{evans02KNNmoments} derives asymptotic limits and convergence rates for
moments of $k$-NN distances, for sampling densities with bounded derivatives
and compact domain. In contrast, we use weaker assumptions to simply prove
bounds on the moments of $k$-NN distances. Importantly, whereas the results of
\citet{evans02KNNmoments} apply only to non-negative moments (i.e.,
$\E \left[ |X|^\alpha \right]$ with $\alpha \geq 0$), our results also hold for
certain negative moments, which is crucial for our bounds on the variance of
the KL estimator.

\subsection{Other Approaches to Estimating Information Theoretic Functionals}
{\bf Analysis of convergence rates:}
For densities over $\R^D$ satisfying a H\"older smoothness condition
parametrized by $\beta \in (0, \infty)$, the minimax rate for estimating
entropy has been known since \citet{birge95estimation} to be
$O \left( n^{-\min \left\{ \frac{8\beta}{4\beta + D}, 1 \right\}} \right)$ in
mean squared error, where $n$ is the sample size.

Quite recently, there has been much work on analyzing new estimators for
entropy, mutual information, divergences, and other functionals of densities.
Most of this work has been along one of three approaches. One series of papers
\citep{liu12exponential,singh14divergence,singh14densityfuncs} studied
boundary-corrected plug-in approach based on under-smoothed kernel density
estimation. This approach has strong finite sample guarantees, but requires
prior knowledge of the support of the density and can necessitate
computationally demanding numerical integration. A second approach
\citep{krishnamurthy14divergences,kandasamy15vonMises} uses von Mises expansion
to correct the bias of optimally smoothed density estimates. This approach
shares the difficulties of the previous approach, but is statistically more
efficient. Finally, a long line of work
\citep{perez08estimation,pal10estimation,sricharan12ensemble,sricharan10confidence,moon14ensemble}
has studied entropy estimation based on continuum limits of certain properties
of graphs (including $k$-NN graphs, spanning trees, and other sample-based
graphs).

Most of these estimators achieve rates of
$O \left( n^{-\min \left\{ \frac{2\beta}{\beta + D}, 1 \right\}} \right)$
or $O \left( n^{-\min \left\{ \frac{4\beta}{2\beta + D}, 1 \right\}} \right)$.
Only the von Mises approach of \citet{krishnamurthy14divergences} is
known to achieve the minimax rate for general $\beta$ and $D$, but due to its
high computational demand ($O(2^D n^3)$), the authors suggest the use
of other statistically less efficient estimators for moderately sized datasets.
In this paper, we prove that, for $\beta \in (0, 2]$, the KL estimator
converges at the rate
$O \left( n^{-\min \left\{ \frac{4\beta}{2\beta + D}, 1 \right\}} \right)$.
It is also worth noting the relative computational efficiency of the KL
estimator ($O \left( D n^2 \right)$, or $O \left( 2^D n \log n \right)$ using
$k$-d trees for small $D$).

{\bf Boundedness of the density:}
For all of the above approaches, theoretical finite-sample results known so far
assume that the sampling density is lower and upper bounded by positive
constants. This also excludes most distributions with unbounded support, and
hence, many distributions of practical relevance. A distinctive feature of our
results is that they hold for a variety of densities that approach $0$ and
$\infty$ on their domain, which may be unbounded. Our bias bounds apply, for
example, to densities that decay exponentially, such as Gaussian distributions.
To our knowledge, the only previous results that apply to unbounded densities
are those of \citet{tsybakov96rootn}, who show $\sqrt{n}$-consistency of a
truncated modification of the KL estimate for a class of functions with
exponentially decaying tails. In fact, components of our analysis are inspired
by \citet{tsybakov96rootn}, and some of our assumptions are closely
related. Their analysis only applies to the case $\beta = 2$ and $D = 1$, for
which our results also imply $\sqrt{n}$-consistency, so our results can be seen
in some respects as a generalization of this work.

\section{Setup and Assumptions}
\label{sec:setting}

While most prior work on $k$-NN estimators has been restricted to $\R^D$, we
present our results in a more general setting. This includes, for example,
Riemannian manifolds embedded in higher dimensional spaces, in which case we
note that our results depend on the \emph{intrinsic}, rather than
\emph{extrinsic}, dimension. Such data can be better behaved in their native
space than when embedded in a lower dimensional Euclidean space (e.g., working
directly on the unit circle avoids boundary bias caused by mapping data to the
interval $[0, 2\pi]$).

\begin{definition}
{\bf (Metric Measure Space):}
A quadruple $(\XX, d, \Sigma, \mu)$ is called a \emph{metric measure space} if
$\XX$ is a set, $d : \XX \times \XX \to [0, \infty)$ is a metric on $\XX$,
$\Sigma$ is a $\sigma$-algebra on $\X$ containing the Borel $\sigma$-algebra
induced by $d$, and $\mu : \Sigma \to [0, \infty]$ is a $\sigma$-finite measure
on the measurable space $(\XX, \Sigma)$.
\label{def:metric_meas_space}
\end{definition}

\begin{definition}
{\bf (Dimension):}
A metric measure space $(\XX, d, \Sigma, \mu)$ is said to have \emph{dimension}
$D \in [0, \infty)$ if there exist constants $c_D, \rho > 0$ such that,
$\forall r \in [0, \rho]$, $x \in \X$, $\mu(B(x, r)) = c_D r^D$.
\footnote{Here and in what follows, $B(x, r) := \{y \in \X : d(x, y) < r\}$
denotes the open ball of radius $r$ centered at $x$.}

\label{def:dim}
\end{definition}

\begin{definition}
{\bf (Full Dimension):}
Given a metric measure space $(\XX, d, \Sigma, \mu)$ of dimension $D$, a
measure $P$ on $(\XX, \Sigma)$ is said to have \emph{full dimension} on a set
$\X \subseteq \XX$ if there exist functions
$\gamma_*, \gamma^* : \X \to (0, \infty)$ such that, for all $r \in [0,\rho]$
and $\mu$-almost all $x \in \X$,
\[\gamma_*(x)r^D \leq P(B(x, r)) \leq \gamma^*(x) r^D.\]
\label{def:full_dim}
\end{definition}

\begin{remark}
If $\XX = \R^D$, $d$ is the Euclidean metric, and $\mu$
is the Lebesgue measure, then the dimension of the metric measure space is $D$.
However, if $\XX$ is a lower dimensional subspace of $\R^D$, then the dimension
may be less than $D$. For example, if
$\XX = \mathbb{S}_{D - 1} := \{x \in \R^D : \|x\|_2 = 1\}$), $d$ is the
geodesic distance on $\mathbb{S}_{D - 1}$, and $\mu$ is the
$(D - 1)$-dimensional surface measure, then the dimension is $D - 1$.
\end{remark}

\begin{remark}
In previous work on $k$-NN statistics
\citep{evans02KNNmoments, biau15EntropyKNN} and estimation of information
theoretic functionals
\citep{sricharan10confidence,krishnamurthy14divergences,singh14divergence,moon14ensemble},
it has been common to make the assumption that the sampling distribution has
full dimension with \emph{constant} $\gamma_*$ and $\gamma^*$ (or, equivalently,
that the density is lower and upper bounded by positive constants). This
excludes distributions with densities approaching $0$ or $\infty$ on their
domain, and hence also densities with unbounded support. By letting $\gamma_*$
and $\gamma^*$ be functions, our results extend to unbounded densities that
instead satisfy certain tail bounds.
\end{remark}

In order to ensure that entropy is well defined, we assume that $P$ is a
probability measure absolutely continuous with respect to $\mu$, and that its
probability density function $p : \X \to [0, \infty)$ satisfies
\footnote{See \cite{baccetti13infiniteEntropy} for discussion of sufficient
conditions for $H(p) < \infty$.}
\begin{equation}
H(p)
  := \E_{X \sim P} \left[ \log p(X) \right]
  = \int_\X p(x) \log p(x) \, d\mu(x) \in \R.
\label{eq:entropy}
\end{equation}
Finally, we assume we have $n + 1$ samples
$X,X_1,...,X_n$ drawn IID from $P$. We would like to use these samples to
estimate the entropy $H(p)$ as defined in Equation (\ref{eq:entropy}).

Our analysis and methods relate to the $k$-nearest neighbor distance $\e_k(x)$,
defined for any $x \in \X$ by $\e_k(x) = d(x, X_i)$, where $X_i$ is the
$k^{th}$-nearest neighbor of $x$ in the set $\{X_1,...,X_n\}$. Note that, since
the definition of dimension used precludes the existence of atoms (i.e.,
for all $x \in \X$, $p(x) = \mu(\{x\}) = 0$), $\e_k(x) > 0$, $\mu$-almost
everywhere. This is important, since we will study $\log \e_k(x)$.

Initially (i.e., in Sections \ref{sec:concentration} and \ref{sec:KNN_stats}),
we will study $\log \e_k(x)$ with fixed $x \in \X$, for which we will derive
bounds in terms of $\gamma_*(x)$ and $\gamma^*(x)$. When we apply these results
to analyze the KL estimator in Section \ref{sec:bias_bound} and
\ref{sec:variance_bound}, we will need to take expectations such as
$\E \left[ \log \e_k(X) \right]$ (for which we reserve the extra sample $X$),
leading to `tail bounds' on $p$ in terms of the functions $\gamma_*$ and
$\gamma^*$.

\section{Concentration of $k$-NN Distances}
\label{sec:concentration}
We begin with a consequence of the multiplicative Chernoff bound, asserting
a sort of concentration of the distance of any point in $\X$ from its
$k^{th}$-nearest neighbor in $\{X_1,\dots,X_n\}$. Since the results of this
section are concerned with fixed $x \in \X$, for notational simplicity, we
suppress the dependence of $\gamma_*$ and $\gamma^*$ on $x$.

\begin{lemma}
Let $(\XX, d, \Sigma, \mu)$ be a metric measure space of dimension $D$. Suppose
$P$ is an absolutely continuous probability measure with full dimension on
$\X \subseteq \XX$ and density function $p : \X \to [0, \infty)$.
For $x \in \X$, if
$r \in \left[\left( \frac{k}{\gamma_* n} \right)^{1/D}, \rho \right]$, then
\[\pr \left[ \e_k(x) > r \right]
  \leq e^{-\gamma_* r^D n} \left( e\frac{\gamma_* r^D n}{k} \right)^k.\]
and, if
$r \in \left[ 0,
              \min\left\{\left( \frac{k}{\gamma^* n} \right)^{1/D},
                         \rho \right\} \right]$, then
\[\pr \left[ \e_k(x) \leq r \right]
  \leq \left( \frac{e \gamma^* r^D n}{k} \right)^{k\gamma_*/\gamma^*}.\]
\label{lemma:KNN_concentration}
\end{lemma}

\section{Bounds on Expectations of KNN Statistics}
\label{sec:KNN_stats}

Here, we use the concentration bounds of Section \ref{sec:concentration} to
bound expectations of functions of $k$-nearest neighbor distances. Specifically,
we give a simple formula for deriving bounds that applies to many functions of
interest, including logarithms and (positive and negative) moments. As in the
previous section, the results apply to a fixed $x \in \X$, and we continue to
suppress the dependence of $\gamma_*$ and $\gamma^*$ on $x$.

\begin{theorem}
Let $(\X, d, \Sigma, \mu)$ be a metric measure space of dimension $D$. Suppose
$P$ is an absolutely continuous probability measure with full dimension and
density function $p : \X \to [0, \infty)$ that satisfies the tail condition
\footnote{Since $f$ need not be surjective, we use the generalized inverse
$f\inv : \R \to [0, \infty]$ defined by
$f\inv(\e) := \inf \{x \in (0, \infty) : f(x) \geq \e\}$.}
\begin{equation}
\E_{X \sim P}
  \left[ \int_\rho^\infty \left[ 1 - P(B(X, f\inv(r))) \right]^n  \right]
  \leq \frac{C_T}{n}
\label{eq:tail_condition}
\end{equation}
for some constant $C_T > 0$. Suppose $f : (0, \infty) \to \R$ is continuously
differentiable, with $f' > 0$. Fix $x \in \X$. Then, we have the upper bound
\begin{align}
\label{ineq:KNN_functional_upper}
& \E \left[ f_+(\e_k(x)) \right]
  \leq f_+\left( \left( \frac{k}{\gamma_* n} \right)^{\frac{1}{D}} \right)
  + \frac{C_T}{n} \\
& + \frac{(e/k)^k}{D(n\gamma_*)^{\frac{1}{D}}}
    \int_k^\infty
        e^{-y} y^{k + \frac{1}{D} - 1}
        f' \left( \left( \frac{y}{n \gamma_*} \right)^{\frac{1}{D}} \right) 
    \, dy
\notag
\end{align}
and the lower bound
\begin{align}
\notag
& \E \left[ f_-(\e_k(x)) \right]
  \leq f_-\left( \left( \frac{k}{\gamma^* n} \right)^{1/D} \right)
  + \frac{C_T}{n} \\
& + \left( \frac{e n \gamma^*}{k} \right)^{\frac{k\gamma_*}{\gamma^*}}
    \int_0^{\left( \frac{k}{\gamma^* n} \right)^{\frac{1}{D}}}
        y^{Dk\gamma_*/\gamma^*} f'(y)
    \, dy
\label{ineq:KNN_functional_lower}
\end{align}
($f_+(x) = \max\{0, f(x)\}$ and $f_-(x) = -\min\{0, f(x)\}$ denote the
positive and negative parts of $f$, respectively).
\label{thm:KNN_functional}
\end{theorem}

\begin{remark}
If $f : (0, \infty) \to \R$ is continuously differentiable with
$f' < 0$, we can apply Theorem \ref{thm:KNN_functional} to $-f$. Also, similar
techniques can be used to prove analogous lower bounds (i.e., lower bounds on
the positive part and upper bounds on the negative part).
\end{remark}
\begin{remark}
The tail condition
(\ref{eq:tail_condition}) is difficult to validate directly for many
distributions. Clearly, it is satisfied when the support of $p$ is
bounded. However, \cite{tsybakov96rootn} show that, for the functions $f$ we
are interested in (i.e., logarithms and power functions), when $\X = \R^D$, $d$
is the Euclidean metric, and $\mu$ is the Lebesgue measure,
(\ref{eq:tail_condition}) is also satisfied by upper-bounded densities with
exponentially decreasing tails. More precisely, that is when there exist
$a,b,\alpha,\delta > 0$ and $\beta > 1$ such that, whenever $\|x\|_2 > \d$,
\[a e^{-\alpha \|x\|^\beta}
  \leq p(x)
  \leq b e^{-\alpha \|x\|^\beta},\]
which permits, for example, Gaussian distributions. It should be noted that the
constant $C_T$ depends only on the metric measure space, the distribution $P$,
and the function $f$, and, in particular, not on $k$.
\end{remark}

\subsection{Applications of Theorem \ref{thm:KNN_functional}}
We can apply Theorem \ref{thm:KNN_functional} to several functions $f$ of
interest. Here, we demonstrate the cases $f(x) = \log x$ and $f(x) = x^\alpha$
for certain $\alpha$, as we will use these bounds when analyzing the KL
estimator.

When $f(x) = \log(x)$, (\ref{ineq:KNN_functional_upper}) gives
\begin{align}
\notag
\E \left[ \log_+(\e_k(x)) \right]
& \leq \frac{1}{D} \log_+ \left( \frac{k}{\gamma_* n} \right)
  + \left( \frac{e}{k} \right)^k \frac{\Gamma(k, k)}{D} \\
& \leq \frac{1}{D}
       \left( 1 + \log_+ \left( \frac{k}{\gamma_* n} \right) \right)
\label{ineq:pos_log_stat}
\end{align}
(where $\Gamma(s, x) := \int_x^\infty t^{s - 1} e^{-t} \, dt$ denotes the upper
incomplete Gamma function, and we used the bound
$\Gamma(s, x) \leq x^{s - 1}e^{-x}$), and (\ref{ineq:KNN_functional_lower})
gives
\begin{align}
\E \left[ \log_-(\e_k(x)) \right]
& \leq \frac{1}{D} \log_-\left( \frac{k}{\gamma^* n} \right)
  + C_1,
\label{ineq:neg_log_stat}
\end{align}
for $C_1 = \frac{\gamma^* e^{k\gamma_*/\gamma^*}}{Dk\gamma_*}$.
For $\alpha > 0$, $f(x) = x^\alpha$, (\ref{ineq:KNN_functional_upper})
gives
\begin{align}
\notag
\E \left[ \e_k^\alpha(x) \right]
& \leq \left( \frac{k}{\gamma_* n} \right)^{\frac{\alpha}{D}}
    + \left( \frac{e}{k} \right)^k
      \frac{\alpha\Gamma\left( k + \alpha/D, k \right)}
           {D(n\gamma_*)^{\alpha/D}} \\
& \leq C_2
       \left( \frac{k}{\gamma_* n} \right)^{\frac{\alpha}{D}},
\label{ineq:pos_moment_stat}
\end{align} 
where $C_2 = 1 + 2\frac{\alpha}{D}$. For any $\alpha \in [-Dk\gamma_*/\gamma^*, 0]$, when $f(x) = -x^\alpha$,
(\ref{ineq:KNN_functional_lower}) gives
\begin{align}
\E \left[ \e_k^\alpha(x) \right]
& \leq C_3
    \left( \frac{k}{\gamma^* n} \right)^{\frac{\alpha}{D}},
\label{ineq:neg_moment_stat}
\end{align}
where $C_3 = 1 + \frac{\alpha \gamma^* e^{k\gamma_*/\gamma^*}}{Dk\gamma_* + \alpha\gamma^*}$.

\section{The KL Estimator for Entropy}
\label{sec:KL_est}
Recall that, for a random variable $X$ sampled from a probability density $p$
with respect to a base measure $\mu$, the Shannon entropy is defined as
\[H(X) = -\int_\X p(x) \log p(x) \, dx.\]
As discussed in Section \ref{sec:intro}, many applications call for estimate of
$H(X)$ given $n$ IID samples $X_1,\dots,X_n \sim p$.
For a positive integer $k$, the KL estimator is typically written as
\[\hat H_k(X)
  = \psi(n) - \psi(k) + \log c_D + \frac{D}{n} \sum_{i = 1}^n \log \e_k(X_i),\]
where $\psi : \N \to \R$ denotes the digamma function. The motivating insight
is the observation that, independent of the sampling distribution,
\footnote{See \cite{Kraskov04estimating} for a concise proof of this fact.}
\[\E \left[ \log P(B(X_i, \e_k(X_i))) \right] = \psi(k) - \psi(n),\]
Hence,
\begin{align*}
& \E \left[ \hat H_k(X) \right] \\
& = \E \left[
        -\log P(B(X_i, \e_k(X_i)))
  + \log c_D + \frac{D}{n} \sum_{i = 1}^n \log \e_k(X_i) \right]  \\
& = -\E \left[ \frac{1}{n} \sum_{i = 1}^n \log \left(
          \frac{P(B(x_i, \e_k(X_i)))}{c_D\e_k^D(X_i)} \right) \right] \\
& = -\E \left[ \frac{1}{n} \sum_{i = 1}^n \log p_{\e_k(i)}(X_i) \right]
  = -\E \left[ \log p_{\e_k(X_1)}(X_1) \right],
\end{align*}
where, for any $x \in \X$, $\e > 0$,
\[p_\e(x)
  = \frac{1}{c_D \e^D} \int_{B(x, \e)} p(y) \, d\mu(y)
  = \frac{P(B(x, \e))}{c_D \e^D}\]
denotes the local average of $p$ in a ball of radius $\e$ around $x$. Since
$p_\e$ is a smoothed approximation of $p$ (with smoothness increasing with
$\e$), the KL estimate can be intuitively thought of as a plug-in estimator for
$H(X)$, using a density estimate with an adaptive smoothing parameter.

In the next two sections, we utilize the bounds derived in Section
\ref{sec:KNN_stats} to bound the bias and variance of the KL estimator. We note
that, for densities in the $\beta$-H\"older smoothness class
($\beta \in (0, 2]$), our results imply a mean-squared error of
$O(n^{-2\beta/D})$ when $\beta < D/2$ and $O(n\inv)$ when $\beta \geq D/2$.

\section{Bias Bound}
\label{sec:bias_bound}

In this section, we prove bounds on the bias of the KL estimator, first in a
relatively general setting, and then, as a corollary, in a more specific but
better understood setting.

\begin{theorem}
Suppose $(\XX, d, \Sigma, \mu)$ and $P$ satisfy the conditions of Theorem
\ref{thm:KNN_functional}, and there exist $C, \beta \in (0, \infty)$ with
\[\sup_{x \in \X} \left| p(x) - p_\e(x) \right| \leq C_\beta \e^\beta,\]
and suppose $p$ satisfies a `tail bound'
\begin{equation}
\Gamma_B
  := \E_{X \sim P} \left[ \left( \gamma_*(X) \right)^{-\frac{\beta + D}{D}} \right]
  < \infty.
\label{ineq:tail_cond}
\end{equation}
Then,
\[\left| \E \left[ H(X) - \hat H_k(X) \right] \right|
  \leq C_B \left( \frac{k}{n} \right)^{\frac{\beta}{D}},\]
where $C_B = (1 + c_D) C_2 C_\beta \Gamma_B$.
\label{thm:gen_bias_bound}
\end{theorem}

We now show that the conditions of Theorem \ref{thm:gen_bias_bound} are
satisfied by densities in the commonly used nonparametric class of
$\beta$-H\"older continuous densities on $\R^D$.

\begin{definition}
Given a constant $\beta > 0$ and an open set $\X \subseteq \R^D$, a function
$f : \X \to \R$ is called \emph{$\beta$-H\"older continuous} if $f$ is $\ell$
times differentiable and there exists $L > 0$ such that, for any multi-index
$\alpha \in \N^D$ with $|\alpha| < \beta$,
\[\sup_{x \neq y \in \X}
      \frac{|D^\alpha f(x) - D^\alpha f(y)|}{\|x - y\|^{\beta - \ell}}
  \leq L,\]
where $\ell := \lfloor \beta \rfloor$ is the greatest integer \emph{strictly}
less than $\beta$.
\end{definition}

\begin{definition}
Given an open set $\X \subseteq \R^D$ and a function $f : \X \to \R$, $f$ is
said to \emph{vanish on the boundary $\partial \X$} of $\X$ if, for any sequence
$\{x_i\}_{i = 1}^\infty$ in $\X$ with
$\inf_{x' \in \partial \X} \|x - x'\|_2 \to 0$ as $i \to \infty$, $f(x) \to 0$
as $i \to \infty$. Here,
\[\partial \X
  := \{x \in \R^D : \forall \d > 0,
                    B(x, \delta) \not\subseteq \X
                    \mbox{ and }
                    B(x, \delta) \not\subseteq \X^c\},\]
denotes the boundary of $\X$.
\end{definition}

\begin{corollary}
Consider the metric measure space $(\R^D, d, \Sigma, \mu)$, where $d$ is Euclidean
and $\mu$ is the Lebesgue measure. Let $P$ be an absolute continuous probability
measure with full dimension and density $p$ supported on an open set
$\X \subseteq \R^D$. Suppose $p$ satisfies (\ref{ineq:tail_cond}) and the
conditions of Theorem \ref{thm:KNN_functional} and is $\beta$-H\"older
continuous ($\beta \in (0, 2]$) with constant $L$. Assume $p$ vanishes on
$\partial \X$. If $\beta > 1$, assume $\|\nabla p\|_2$ vanishes on
$\partial \X$. Then,
\[\left| \E \left[ \hat H_k(X) - H(X) \right] \right|
  \leq C_H \left( \frac{n}{k} \right)^{-\frac{\beta}{D}},\]
where $C_H = (1 + c_D) C_2 \Gamma \frac{LD}{D + \beta}$.
\label{corr:Holder_bias_bound}
\end{corollary}

\begin{remark}
The assumption that $p$ (and perhaps $\|\nabla p\|$) vanish on the boundary of
$\X$ can be thought of as ensuring that the trivial continuation
$q : \R^D \to [0, \infty)$
\[q(x)
  = \left\{
      \begin{array}{ll}
        p(x) & x \in \X \\
        0 & x \in \R^D \sminus \X
      \end{array}
    \right.\]
of $p$ to $\R^D$ is $\beta$-H\"older continuous. This reduces boundary bias,
for which the KL estimator does not correct.
\footnote{Several estimators controlling for boundary bias have been proposed
(e.g., \citet{sricharan10confidence} give a modified $k$-NN estimator that
accomplishes this \emph{without} prior knowledge of $\X$.}
\end{remark}

\section{Variance Bound}
\label{sec:variance_bound}

We first use the bounds proven in Section \ref{sec:KNN_stats} to prove uniform
(in $n$) bounds on the moments of $\E \left[ \log \e_k(X) \right]$. We the
for any fixed $x \in \X$, although $\log \e_k(x) \to -\infty$ almost surely as
$n \to \infty$, $\Var \left[ \log \e_k(x) \right]$, and indeed all higher
central moments of $\log \e_k(x)$, are bounded, uniformly in $n$. In fact,
there exist exponential bounds, independent of $n$, on the
density of $\log \e_k(x) - \E \left[ \log \e_k(x) \right]$.

\subsection{Moment Bounds on Logarithmic $k$-NN distances}

\begin{lemma}
Suppose $(\XX, d, \Sigma, \mu)$ and $P$ satisfy the conditions of Theorem
\ref{thm:KNN_functional}. Suppose also that
$\Gamma_0 := \sup_{x \in \X} \frac{\gamma^*(x)}{\gamma_*(x)} < \infty$. Let
$\lambda \in \left(0, \frac{Dk}{\Gamma_0} \right)$ and assume the following
expectations are finite:
\begin{equation}
\Gamma := \E_{X \sim P} \left[ \frac{\gamma^*(X)}{\gamma_*(X)} \right] < \infty.
\label{const:gamma_ratio}
\end{equation}
\begin{equation}
\Gamma_*(\lambda) := \E_{X \sim P} \left[ \left( \gamma_*(X) \right)^{-\lambda/D} \right] < \infty.
\label{const:gamma_sub_pow}
\end{equation}
\begin{equation}
\Gamma^*(\lambda) := \E_{X \sim P} \left[ \left( \gamma^*(X) \right)^{\lambda/D} \right] < \infty.
\label{const:gamma_sup_pow}
\end{equation}

Then, for any integer $\ell > 1$, the $\ell^{th}$ central moment
\[M_\ell := \E \left[ \left( \log \e_k(X) - \E \left[ \log \e_k(X) \right] \right)^\ell \right]\]
satisfies
\begin{equation}
M_\ell \leq C_M \ell!/\lambda^\ell,
\label{ineq:moment_bound}
\end{equation}
where $C_M > 0$ is a constant independent of $n$, $\ell$, and $\lambda$.
\label{lemma:log_moment_bound}
\end{lemma}

\begin{remark}
The conditions (\ref{const:gamma_ratio}), (\ref{const:gamma_sub_pow}), and
(\ref{const:gamma_sup_pow}) are mild. For example, when $\X = \R^D$, $d$ is the
Euclidean metric, and $\mu$ is the Lebesgue measure, it suffices that $p$ is
Lipschitz continuous
\footnote{Significantly milder conditions than Lipschitz continuity suffice,
but are difficult to state here due to space limitations.} and there exist
$c, r > 0, p > \frac{D^2}{D - \alpha}$ such that $p(x) \leq c\|x\|^{-p}$
whenever $\|x\|_2 > r$. The condition $\Gamma_0 < \infty$ is more prohibitive,
but still permits many (possibly unbounded) distributions of interest.
\end{remark}

\begin{remark}
If the terms $\log \e_k(X_i)$ were independent, a Bernstein inequality,
together with the moment bound (\ref{ineq:moment_bound}) would imply a
sub-Gaussian concentration bound on the KL estimator about its expectation.
This may follow from one of several more refined concentration results relaxing
the independence assumption that have been proposed.
\end{remark}

\subsection{Bound on the Variance of the KL Estimate}
Bounds on the variance of the KL estimator now follow from the law of large
numbers in \citet{evans08SLLNforKNN} (itself an application of the Efron-Stein
inequality to $k$-NN statistics).
\begin{theorem}
Suppose $(\XX, d, \Sigma, \mu)$ and $P$ satisfy the conditions of Lemma
\ref{lemma:log_moment_bound}, and that that there exists a constant $N_k \in \N$
such that, for any finite $F \subseteq \X$, any $x \in F$ can be among the
$k$-NN of at most $N_k$ other points in that set. Then,
$\hat H_k(X) \to \E \left[ \hat H_k(X) \right]$ almost surely (as
$n \to \infty$), and, for $n \geq 16 k$ and $M_4$ satisfying
(\ref{ineq:moment_bound}).
\[\Var \left[ \hat H_k(X) \right]
  \leq \frac{5(3 + kN_k)(3 + 64k)M_4}{n}
  \in O \left( \frac{1}{nk} \right),\]
\label{thm:variance_bound}
\end{theorem}

\begin{remark}
$N_k$ depends only on $k$ and the geometry of the metric space $(\X, d)$. For
example, Corollary A.2 of \citet{evans08SLLNforKNN} shows that, when
$\X = \R^D$ and $d$ is the Euclidean metric, then $N_k \leq k K(D)$, where
$K(D)$ is the kissing number of $\R^d$.
\end{remark}

\section{Bounds on the Mean Squared Error}
The bias and variance bounds (Theorems \ref{thm:gen_bias_bound} and
\ref{thm:variance_bound}) imply a bound on the mean squared error of the KL
estimator:
\begin{corollary}
Suppose $p$
\begin{enumerate}
\item
is $\beta$-H\"older continuous with $\beta \in (0, 2]$.
\item
vanishes on $\partial \X$. If $\beta > 1$, then also suppose
$\|\nabla p\|_2$ vanishes on $\partial \X$.
\item
\end{enumerate}
[TODO: Other assumptions.]
satisfies the assumptions of Theorems \ref{thm:gen_bias_bound} and
\ref{thm:variance_bound}. Then,
\begin{equation}
\E \left[ \left( \hat H_k(X) - H(X) \right)^2 \right]
  \leq C_B^2 \left( \frac{k}{n} \right)^{2\beta/D} + \frac{C_V}{nk}.
\label{ineq:MSE_bound_general_k}
\end{equation}
If we let $k$ scale as $k \asymp n^{\max \left\{ 0, \frac{2 \beta - D}{2 \beta + D} \right\}}$
this gives an overall convergence rate of
\begin{equation}
\E \left[ \left( \hat H_k(X) - H(X) \right)^2 \right]
  \leq C_B^2 \left( \frac{k}{n} \right)^{2\beta/D} + \frac{C_V}{nk}.
\label{ineq:MSE_bound_optimal_k}
\end{equation}
\label{corr:MSE_bound}
\end{corollary}

\section{Conclusions and Future Work}

This paper derives finite sample bounds on the bias and variance of the KL
estimator under general conditions, including for certain classes of unbounded
distributions. As intermediate results, we proved concentration inequalities for
$k$-NN distances and bounds on the expectations of statistics of $k$-NN
distances. We hope these results and methods may lead to convergence rates for
the widely used KSG mutual information estimator, or to generalize convergence
rates for other estimators of entropy and related functionals to unbounded
distributions.

\section*{Acknowledgements} 
This material is based upon work supported by a National Science Foundation
Graduate Research Fellowship to the first author under Grant No. DGE-1252522.

\bibliography{biblio}
\bibliographystyle{icml2016}

\end{document}